\documentclass[a4paper,leqno,12pt]{article}

\usepackage{amsmath}
\usepackage{amssymb}
\usepackage{amsthm}

\begin{document}

\bibliographystyle{amsalpha}
\newtheorem{Assumption}{Assumption}[section]
\newtheorem{Theorem}{Theorem}[section]
\newtheorem{Lemma}{Lemma}[section]
\newtheorem{Remark}{Remark}[section]
\newtheorem{Corollary}{Corollary}[section]
\newtheorem{Conjecture}{Conjecture}[section]
\newtheorem{Proposition}{Proposition}[section]
\newtheorem{Example}{Example}[section]
\newtheorem{Definition}{Definition}[section]
\newtheorem{Problem}{Problem}[section]
\newtheorem{Subsection}{Subsection}[section]
\newtheorem{Proof}{Proof}[section]

\title{Studies on the Lorenz model \\}

\author{By\\
Yusuke Sasano \\
(The University of Tokyo, Japan)}
\maketitle


\begin{abstract} We study the Lorenz model from the viewpoint of its accessible singularities and local index.


\end{abstract}

\section{Introduction}

The Lorenz model
\begin{equation}\label{l1}
  \left\{
  \begin{aligned}
   \frac{dx}{dt} &=y-\sigma \varepsilon x,\\
   \frac{dy}{dt} &=-xz+x-\varepsilon y,\\
   \frac{dz}{dt} &=xy-\varepsilon b z\\
   \end{aligned}
  \right. 
\end{equation}
has been studied by Tabor and Weiss (1981) in detail \cite{15}. In this paper, we study the phase space of \eqref{l1} from the viewpoint of its accessible singularities and local index.

As the conditions which each accessible singular point can be resolved, we obtain the following
  \begin{align}
  \begin{split}
    & \varepsilon(b-1)(b-2\sigma)(b+3\sigma-1)=0,\\
    & (b-1)(b-3\sigma+1)=0,\\
    & \{b^2-5b-2-3(b-2)\sigma \} \{\varepsilon^2(b-1)(7b-15\sigma+2)-9\}=0,\\
    & \varepsilon(b-2\sigma)(b+3\sigma-1)=0.\\
    \end{split}
   \end{align}
These equations can be solved as follows:
\begin{equation*}
\{(\sigma,\varepsilon,b)=\left(\frac{1}{3},\varepsilon,0 \right),(1,-3,2),(1,3,2),(2,0,1)\}.
\end{equation*}
In each case, we study its first integrals, general solutions and phase space.

\section{Accessible singularities }

Let us review the notion of accessible singularity. Let $B$ be a connected open domain in $\Bbb C$ and $\pi : {\mathcal W} \longrightarrow B$ a smooth proper holomorphic map. We assume that ${\mathcal H} \subset {\mathcal W}$ is a normal crossing divisor which is flat over $B$. Let us consider a rational vector field $\tilde v$ on $\mathcal W$ satisfying the condition
\begin{equation*}
\tilde v \in H^0({\mathcal W},\Theta_{\mathcal W}(-\log{\mathcal H})({\mathcal H})).
\end{equation*}
Fixing $t_0 \in B$ and $P \in {\mathcal W}_{t_0}$, we can take a local coordinate system $(x_1,\ldots ,x_n)$ of ${\mathcal W}_{t_0}$ centered at $P$ such that ${\mathcal H}_{\rm smooth \rm}$ can be defined by the local equation $x_1=0$.
Since $\tilde v \in H^0({\mathcal W},\Theta_{\mathcal W}(-\log{\mathcal H})({\mathcal H}))$, we can write down the vector field $\tilde v$ near $P=(0,\ldots ,0,t_0)$ as follows:
\begin{equation*}
\tilde v= \frac{\partial}{\partial t}+a_1 
\frac{\partial}{\partial x_1}+\frac{a_2}{x_1} 
\frac{\partial}{\partial x_2}+\cdots+\frac{a_n}{x_1} 
\frac{\partial}{\partial x_n}.
\end{equation*}
This vector field defines the following system of differential equations
\begin{equation}\label{39}
  \left\{
  \begin{aligned}
   \frac{dx_1}{dt} &=a_1(x_1,x_2,....,x_n,t),\\
   \frac{dx_2}{dt} &=\frac{a_2(x_1,x_2,....,x_n,t)}{x_1},\\
   .\\
   .\\
   .\\
   \frac{dx_n}{dt} &=\frac{a_n(x_1,x_2,....,x_n,t)}{x_1}.
   \end{aligned}
  \right. 
\end{equation}
Here $a_i(x_1,\ldots,x_n,t), \ i=1,2,\dots ,n,$ are holomorphic functions defined near $P=(0,\dots ,0,t_0).$

\begin{Definition}\label{Def}
With the above notation, assume that the rational vector field $\tilde v$ on $\mathcal W$ satisfies the condition
$$
(A) \quad \tilde v \in H^0({\mathcal W},\Theta_{\mathcal W}(-\log{\mathcal H})({\mathcal H})).
$$
We say that $\tilde v$ has an {\it accessible singularity} at $P=(0,\dots ,0,t_0)$ if
$$
x_1=0 \ {\rm and \rm} \ a_i(0,\ldots,0,t_0)=0 \ {\rm for \rm} \ {\rm every \rm} \ i, \ 2 \leq i \leq n.
$$
\end{Definition}

If $P \in {\mathcal H}_{{\rm smooth \rm}}$ is not an accessible singularity, all solutions of the ordinary differential equation passing through $P$ are vertical solutions, that is, the solutions are contained in the fiber ${\mathcal W}_{t_0}$ over $t=t_0$. If $P \in {\mathcal H}_{\rm smooth \rm}$ is an accessible singularity, there may be a solution of \eqref{39} which passes through $P$ and goes into the interior ${\mathcal W}-{\mathcal H}$ of ${\mathcal W}$.

Here we review the notion of {\it local index}. Let $v$ be an algebraic vector field with an accessible singular point $\overrightarrow{p}=(0,\ldots,0)$ and $(x_1,\ldots,x_n)$ be a coordinate system in a neighborhood centered at $\overrightarrow{p}$. Assume that the system associated with $v$ near $\overrightarrow{p}$ can be written as
\begin{align}\label{b}
\begin{split}
\frac{d}{dt}Q\begin{pmatrix}
             x_1 \\
             x_2 \\
             \vdots\\
             x_n
             \end{pmatrix}=\frac{1}{x_1}\{Q\begin{bmatrix}
             a_1 & & & &  \\
             & a_2 & & & \\
             & & \ddots & & \\
             & & & & a_n
             \end{bmatrix}Q^{-1}{\cdot}Q\begin{pmatrix}
             x_1 \\
             x_2 \\
             \vdots\\
             x_n
             \end{pmatrix}+\begin{pmatrix}
             x_1f_1(x_1,x_2,\ldots,x_n,t) \\
             f_2(x_1,x_2,\ldots,x_n,t) \\
             \vdots\\
             f_n(x_1,x_2,\ldots,x_n,t)
             \end{pmatrix}\},\\
              (f_i \in {\mathbb C}(t)[x_1,\ldots,x_n], \ Q \in GL(n,{\mathbb C}(t)), \ a_i \in {\mathbb C}(t))
             \end{split}
             \end{align}
where $f_1$ is a polynomial which vanishes at $\overrightarrow{p}$ and $f_i$, $i=2,3,\ldots,n$ are polynomials of order at least 2 in $x_1,x_2,\ldots,x_n$. We call ordered set of the eigenvalues $(a_1,a_2,\ldots,a_n)$ {\it local index} at $\overrightarrow{p}$.

We remark that we are interested in the case with local index
\begin{equation}
(1,a_2/a_1,\ldots,a_n/a_1) \in {\mathbb Z}^n.
\end{equation}
If each component of $(1,a_2/a_1,\ldots,a_n/a_1)$ has the same sign, we may resolve the accessible singularity by blowing-up finitely many times. However, when different signs appear, we may need to both blow up and blow down.

In order to consider the phase spaces for the system \eqref{l1}, let us take the compactification $[z_0:z_1:z_2:z_3] \in {\mathbb P}^3$ of $(x,y,z) \in {\mathbb C}^3$ with the natural embedding
$$
(x,y,z)=(z_1/z_0,z_2/z_0,z_3/z_0).
$$
Moreover, we denote the boundary divisor in ${\mathbb P}^3$ by $ {\cal H}$. Extend the regular vector field on ${\mathbb C}^3$ to a rational vector field $\tilde v$ on ${\mathbb P}^3$. It is easy to see that ${\mathbb P}^3$ is covered by four copies of ${\mathbb C}^3${\rm : \rm}
\begin{align*}
&U_0={\mathbb C}^3 \ni (x,y,z),\\
&U_j={\mathbb C}^3 \ni (X_j,Y_j,Z_j) \ (j=1,2,3),
\end{align*}
via the following rational transformations
\begin{align*}
& X_1=1/x, \quad Y_1=y/x, \quad Z_1=z/x,\\
& X_2=x/y, \quad Y_2=1/y, \quad Z_2=z/y,\\
& X_3=x/z, \quad Y_3=y/z, \quad Z_3=1/z.
\end{align*}
The following Lemma shows that this rational vector field $\tilde v$ has five accessible singular points on the boundary divisor ${\mathcal H} \subset {\mathbb P}^3$.
\begin{Lemma}
The rational vector field $\tilde v$ has five accessible singular points{\rm : \rm}
\begin{equation}
  \left\{
  \begin{aligned}
   P_1 &=\{(X_1,Y_1,Z_1)|X_1=Y_1=Z_1=0\},\\
   P_2 &=\{(X_2,Y_2,Z_2)|X_2=Y_2=Z_2=0\},\\
   P_3 &=\{(X_3,Y_3,Z_3)|X_3=Y_3=Z_3=0\},\\
   P_4 &=\{(X_2,Y_2,Z_2)|X_2=Y_2=0, \ Z_2=\sqrt{-1}\},\\
   P_5 &=\{(X_2,Y_2,Z_2)|X_2=Y_2=0, \ Z_2=-\sqrt{-1}\}.
   \end{aligned}
  \right. 
\end{equation}
\end{Lemma}

Next let us calculate its local index at the points $P_i \ (i=1,2,3)$.
\begin{center}
\begin{tabular}{|c|c|c|} \hline 
Singular point & Type of local index   \\ \hline 
$P_1$ & $(0,\sqrt{-1},-\sqrt{-1})$  \\ \hline 
$P_2$ & $(0,0,0)$  \\ \hline 
$P_3$ & $(0,0,0)$  \\ \hline 
\end{tabular}
\end{center}
We see that there are no solutions which pass through $P_i, \ (i=1,2,3)$, respectively.

In order to do analysis for the accessible singularities $P_4,P_5$, we need to replace a suitable coordinate system because each point has multiplicity of order 2.

At first, let us do the Painlev\'e test. To find the leading order behaviour of a singularity at $t=t_1$ one sets
\begin{equation*}
  \left\{
  \begin{aligned}
   x & \propto \frac{a}{(t-t_1)^m},\\
   y & \propto \frac{b}{(t-t_1)^n},\\
   z & \propto \frac{c}{(t-t_1)^p},
   \end{aligned}
  \right. 
\end{equation*}
from which it is easily deduced that
\begin{equation*}
m=1, \quad n=2, \quad p=2
\end{equation*}
and
\begin{equation*}
a=\pm 2\sqrt{-1}, \quad n=\mp 2\sqrt{-1}, \quad p=-2.
\end{equation*}
Each order of pole $(m,n,p)$ suggests a suitable coordinate system to do analysis for the accessible singularities $P_4,P_5$, which is explicitly given by
\begin{equation*}
(X,Y,Z)=\left(\frac{1}{x},\frac{y}{x^2},\frac{z}{x^2} \right).
\end{equation*}
In this coordinate, the singular points $P_4,P_5$ are given as follows:
\begin{equation*}
  \left\{
  \begin{aligned}
   P_4 &= \{(X,Y,Z)=\left(0,\frac{\sqrt{-1}}{2},\frac{1}{2} \right)\},\\
   P_5 &= \{(X,Y,Z)=\left(0,-\frac{\sqrt{-1}}{2},\frac{1}{2} \right)\}.
   \end{aligned}
  \right. 
\end{equation*}
Next let us calculate its local index at each point.
\begin{center}
\begin{tabular}{|c|c|c|} \hline 
Singular point & Type of local index   \\ \hline 
$P_4$ & $(-\frac{\sqrt{-1}}{2},-2\sqrt{-1},-\sqrt{-1})$  \\ \hline 
$P_5$ & $(\frac{\sqrt{-1}}{2},2\sqrt{-1},\sqrt{-1})$  \\ \hline 
\end{tabular}
\end{center}

Now, we try to resolve the accessible singular points $P_4,P_5$.

{\bf Step 0}: We take the coordinate system centered at $P_4${\rm : \rm}
$$
p=X, \quad q=Y-\frac{\sqrt{-1}}{2}, \quad r=Z-\frac{1}{2}.
$$

{\bf Step 1}: We make the linear transformation{\rm : \rm}
$$
p^{(1)}=p, \quad q^{(1)}=q-\sqrt{-1} r, \quad r^{(1)}=r.
$$
In this coordinate, the system \eqref{l1} is rewritten as follows:
\begin{align*}
\frac{d}{dt}\begin{pmatrix}
             p^{(1)} \\
             q^{(1)} \\
             r^{(1)} 
             \end{pmatrix}&=\frac{1}{p^{(1)}}\{\begin{pmatrix}
             -\frac{\sqrt{-1}}{2} & 0 & 0 \\
             0 & -2\sqrt{-1} & 0 \\
             0 & 0 & -\sqrt{-1}
             \end{pmatrix}\begin{pmatrix}
             p^{(1)} \\
             q^{(1)} \\
             r^{(1)} 
             \end{pmatrix}+\dots\}.
             \end{align*}
By considering the ratio $\left(1,\frac{-2\sqrt{-1}}{-\frac{\sqrt{-1}}{2}},\frac{-\sqrt{-1}}{-\frac{\sqrt{-1}}{2}} \right)=(1,4,2)$, we obtain the resonances $(4,2)$. This property suggests that we will blow up four times to the direction $q^{(1)}$ and two times to the direction $r^{(1)}$.

{\bf Step 2}: We blow up at the point $P_4=\{(p^{(1)},q^{(1)},r^{(1)})=(0,0,0)\}${\rm : \rm}
$$
p^{(2)}=p^{(1)}, \quad q^{(2)}=\frac{q^{(1)}}{p^{(1)}}, \quad r^{(2)}=\frac{r^{(1)}}{p^{(1)}}.
$$

{\bf Step 3}: We blow up at the point $\{(p^{(2)},q^{(2)},r^{(2)})=(0,\frac{\varepsilon}{3}(b-1),\sqrt{-1} \varepsilon(b-2\sigma)\}${\rm : \rm}
$$
p^{(3)}=p^{(2)}, \quad q^{(3)}=\frac{q^{(2)}-\frac{\varepsilon}{3}(b-1)}{p^{(2)}}, \quad r^{(3)}=\frac{r^{(2)}-\sqrt{-1} \varepsilon(b-2\sigma)}{p^{(2)}}.
$$

{\bf Step 4}: We blow up along the curve $\{(p^{(3)},q^{(3)},r^{(3)})|p^{(3)}=0,\\
q^{(3)}=\frac{\sqrt{-1}}{9}\left(\varepsilon^2(b-1)(7b-15\sigma+2)-9 \right) \}${\rm : \rm}
$$
p^{(4)}=p^{(3)}, \quad q^{(4)}=\frac{q^{(3)}-\frac{\sqrt{-1}}{9}\left(\varepsilon^2(b-1)(7b-15\sigma+2)-9 \right)}{p^{(3)}}, \quad r^{(4)}=r^{(3)}.
$$

{\bf Step 5}: We blow up along the curve $\{(p^{(4)},q^{(4)},r^{(4)})|p^{(4)}=0,\\
q^{(4)}=\frac{4}{3}\varepsilon(b-1)q^{(4)}-\frac{2}{27}\varepsilon(b+2)\left(\varepsilon^2(b-1)(7b-15\sigma+2)-9 \right) \}${\rm : \rm}
$$
u=p^{(4)}, \quad v=\frac{q^{(4)}-\left(\frac{4}{3}\varepsilon(b-1)r^{(4)}-\frac{2}{27}\varepsilon(b+2)\left(\varepsilon^2(b-1)(7b-15\sigma+2)-9 \right) \right)}{p^{(4)}},
$$
$$
w=r^{(4)}.
$$
In this coordinate, the system \eqref{l1} is rewritten as follows:
\begin{equation}\label{A}
  \left\{
  \begin{aligned}
   \frac{du}{dt} &=g_1(u,v,w),\\
   \frac{dv}{dt} &=\frac{4\sqrt{-1}}{9}\frac{\varepsilon^3(b-1)(b-2\sigma)(b+3\sigma-1)}{u^2}-\frac{2\varepsilon^2}{27}\frac{18(b-1)(b-3\sigma+1)w}{u}\\
   &-\frac{2\varepsilon^2}{27}\frac{\{b^2-5b-2-3(b-2)\sigma \}\{\varepsilon^2(b-1)(7b-15\sigma+2)-9\}}{u}+g_2(u,v,w),\\
   \frac{dw}{dt} &=-\frac{\sqrt{-1}}{3}\frac{\varepsilon^2(b-2\sigma)(b+3\sigma-1)}{u}+g_3(u,v,w),
   \end{aligned}
  \right. 
\end{equation}
where $g_i(u,v,w) \in {\mathbb C}[u,v,w] \ (i=1,2,3)$.

Each right-hand side of the system \eqref{A} is a {\it polynomial} if and only if
  \begin{align}\label{B}
  \begin{split}
   &\varepsilon(b-1)(b-2\sigma)(b+3\sigma-1)=0,\\
   &(b-1)(b-3\sigma+1)=0,\\
   &\{b^2-5b-2-3(b-2)\sigma \}\{\varepsilon^2(b-1)(7b-15\sigma+2)-9\}=0,\\
   &\varepsilon(b-2\sigma)(b+3\sigma-1)=0.
   \end{split}
   \end{align}
These equations can be solved as follows:
\begin{equation}\label{C}
\{(\sigma,\varepsilon,b)=\left(\frac{1}{3},\varepsilon,0 \right),(1,-3,2),(1,3,2),(2,0,1)\}.
\end{equation}

\section{The case of $(\sigma,\varepsilon,b)=\left(\frac{1}{3},\varepsilon,0 \right)$}

\begin{equation}\label{2.1}
  \left\{
  \begin{aligned}
   \frac{dx}{dt} &=y-\frac{\varepsilon}{3}x,\\
   \frac{dy}{dt} &=-xz+x-\varepsilon y,\\
   \frac{dz}{dt} &=xy.
   \end{aligned}
  \right. 
\end{equation}

This system is equivalent to the 3-rd order ordinary differential equation:
\begin{equation}
\frac{d^3x}{dt^3}=-\frac{\varepsilon}{3}x^3-x^2\frac{dx}{dt}+\frac{4\varepsilon}{3x}\left(\frac{dx}{dt}\right)^2+\frac{1}{x}\frac{dx}{dt}\frac{d^2x}{dt^2}-\frac{4\varepsilon}{3}\frac{d^2x}{dt^2}.
\end{equation}
This system does not appear in the Chazy polynomial class.

\begin{Theorem}
The phase space ${\mathcal X}$ for the system \eqref{2.1} is obtained by gluing three copies of ${\mathbb C}^3${\rm:\rm}
\begin{center}
${U_j} \cong {\mathbb C}^3 \ni \{(x_j,y_j,z_j)\},  \ \ j=0,1,2$
\end{center}
via the following birational transformations{\rm:\rm}
\begin{align}
\begin{split}
0) \ x_0=&x, \quad y_0=y, \quad z_0=z,\\
1) \ x_1=&\frac{1}{x}, \quad y_1=-\left((y-\frac{\varepsilon}{3}x-\sqrt{-1}z+\frac{\sqrt{-1}(5\varepsilon^2+9)}{9})x+\frac{4\varepsilon}{9}(3z+\varepsilon^2-3)\right)x,\\
&z_1=z-\frac{1}{6}(3x-4\sqrt{-1} \varepsilon)x,\\
2) \ x_2=&\frac{1}{x}, \quad y_2=-\left((y-\frac{\varepsilon}{3}x+\sqrt{-1}z-\frac{\sqrt{-1}(5\varepsilon^2+9)}{9})x+\frac{4\varepsilon}{9}(3z+\varepsilon^2-3)\right)x,\\
&z_2=z-\frac{1}{6}(3x+4\sqrt{-1} \varepsilon)x.
\end{split}
\end{align}
These transition functions satisfy the condition{\rm:\rm}
\begin{equation*}
dx_i \wedge dy_i \wedge dz_i=dx \wedge dy \wedge dz \quad (i=1,2).
\end{equation*}
\end{Theorem}

\begin{Theorem}
Let us consider a system of first order ordinary differential equations in the polynomial class\rm{:\rm}
\begin{equation*}
\frac{dx}{dt}=f_1(x,y,z), \quad \frac{dy}{dt}=f_2(x,y,z), \quad \frac{dz}{dt}=f_3(x,y,z).
\end{equation*}
We assume that

$(A1)$ $deg(f_i)=2$ with respect to $x,y,z$.

$(A2)$ The right-hand side of this system becomes again a polynomial in each coordinate system $(x_i,y_i,z_i) \ (i=1,2)$.

\noindent
Then such a system coincides with the system \eqref{2.1}.
\end{Theorem}

\section{Modified Lorenz model}
In this section, we present 4-parameter family of modified Lorenz model explicitly given by

\begin{equation}\label{m2.1}
  \left\{
  \begin{aligned}
   \frac{dx}{dt} &=-\frac{\varepsilon}{3}x+y+\frac{1}{72}\left(8(9\alpha_1-\varepsilon \alpha_3)+\sqrt{-1}(24\alpha_2+\alpha_3^2-16\varepsilon^2)\right),\\
   \frac{dy}{dt} &=-xz-\frac{1}{72}\left(24\alpha_2+\alpha_3^2-8\varepsilon^2 \right)x-\varepsilon y+\frac{1}{6}(\alpha_3-4\sqrt{-1} \varepsilon)z\\
   &+\frac{1}{432}(24\alpha_2 \alpha_3+\alpha_3^3-432\alpha_1 \varepsilon +64\alpha_3 \varepsilon^2-2\sqrt{-1}\varepsilon(120\alpha_2+5\alpha_3^2-16\varepsilon^2),\\
   \frac{dz}{dt} &=xy+\frac{1}{72}\left(12(6\alpha_1-\varepsilon \alpha_3)+\sqrt{-1}(24\alpha_2+\alpha_3^2)\right)x-\frac{1}{6}(\alpha_3-4\sqrt{-1} \varepsilon)y\\
   &-\frac{1}{432}(\alpha_3-4\sqrt{-1} \varepsilon)(12(6\alpha_1-\alpha_3 \varepsilon)+\sqrt{-1}(24\alpha_2+\alpha_3^2)),
   \end{aligned}
  \right. 
\end{equation}
where $\alpha_i,\varepsilon$ are complex parameters.

\begin{Theorem}
The phase space ${\mathcal X}$ for the system \eqref{m2.1} is obtained by gluing three copies of ${\mathbb C}^3${\rm:\rm}
\begin{center}
${U_j} \cong {\mathbb C}^3 \ni \{(x_j,y_j,z_j)\},  \ \ j=0,1,2$
\end{center}
via the following birational transformations{\rm:\rm}
\begin{align}
\begin{split}
0) \ x_0=&x, \quad y_0=y, \quad z_0=z,\\
1) \ x_1=&\frac{1}{x}, \quad y_1=-\left((y-\frac{\varepsilon}{3}x-\sqrt{-1}z+\alpha_1)x+\frac{4\varepsilon}{9}(3z+\alpha_2)\right)x,\\
&z_1=z-\frac{1}{6}(3x-\alpha_3)x,\\
2) \ x_2=&\frac{1}{x}, \quad y_2=-((y-\frac{\varepsilon}{3}x+\sqrt{-1}z+\frac{1}{36}(36\alpha_1+24\sqrt{-1}\alpha_2+\sqrt{-1}\alpha_3^2-48\sqrt{-1} \varepsilon^2)x\\
&+\frac{4\varepsilon}{9}(3z+\frac{1}{3}(3\alpha_2+2\sqrt{-1}\alpha_3\varepsilon+8\varepsilon^2))x,\\
&z_2=z-\frac{1}{6}(3x-\alpha_3+8\sqrt{-1}\varepsilon)x.
\end{split}
\end{align}
These transition functions satisfy the condition{\rm:\rm}
\begin{equation*}
dx_i \wedge dy_i \wedge dz_i=dx \wedge dy \wedge dz \quad (i=1,2).
\end{equation*}
\end{Theorem}

\begin{Theorem}
Let us consider a system of first order ordinary differential equations in the polynomial class\rm{:\rm}
\begin{equation*}
\frac{dx}{dt}=f_1(x,y,z), \quad \frac{dy}{dt}=f_2(x,y,z), \quad \frac{dz}{dt}=f_3(x,y,z).
\end{equation*}
We assume that

$(A1)$ $deg(f_i)=2$ with respect to $x,y,z$.

$(A2)$ The right-hand side of this system becomes again a polynomial in each coordinate system $(x_i,y_i,z_i) \ (i=1,2)$.

\noindent
Then such a system coincides with the system \eqref{m2.1}.
\end{Theorem}

\section{The case of $(\sigma,\varepsilon,b)=(2,0,1)$}

\begin{equation}\label{3.1}
  \left\{
  \begin{aligned}
   \frac{dx}{dt} &=y,\\
   \frac{dy}{dt} &=-xz+x,\\
   \frac{dz}{dt} &=xy.
   \end{aligned}
  \right. 
\end{equation}

\begin{Proposition}
This system has
\begin{equation}
I:=x^2-2z
\end{equation}
as its first integral.
\end{Proposition}

This system can be solved by reduction to 2nd-order ordinary differential equation
\begin{equation}
  \left\{
  \begin{aligned}
   \frac{dx}{dt} &=y,\\
   \frac{dy}{dt} &=-\frac{1}{2}x^3+(1+\frac{I}{2})x,
   \end{aligned}
  \right. 
\end{equation}
or equivalently
\begin{equation}
\frac{d^2x}{dt^2}=-\frac{1}{2}x^3+(1+\frac{I}{2})x,
\end{equation}
which is a special case of equation Ince-VIII.

\section{The case of $(\sigma,\varepsilon,b)=(1,3,2)$}

\begin{equation}\label{4.1}
  \left\{
  \begin{aligned}
   \frac{dx}{dt} &=y-3x,\\
   \frac{dy}{dt} &=-xz+x-3y,\\
   \frac{dz}{dt} &=xy-6z.
   \end{aligned}
  \right. 
\end{equation}

\begin{Proposition}
This system has
\begin{equation}
I:=e^{6t}(x^2-2z)
\end{equation}
as its first integral.
\end{Proposition}

Owing to its first integral, this system can be reduced to the system
\begin{equation}
  \left\{
  \begin{aligned}
   \frac{dx}{dt} &=y-3x,\\
   \frac{dy}{dt} &=-\frac{1}{2}x^3-3y+\frac{e^{-6t}(I+2e^{6t})}{2}x.\\
   \end{aligned}
  \right. 
\end{equation}
By making the change of variables
\begin{equation}
  \left\{
  \begin{aligned}
   X &=\frac{\sqrt{-1}}{2}x,\\
   Y &=\frac{2x+\sqrt{-1}x^2-2y}{2x},
   \end{aligned}
  \right. 
\end{equation}
we obtain
\begin{equation}
  \left\{
  \begin{aligned}
   \frac{dX}{dt} &=X^2-XY-2X,\\
   \frac{dY}{dt} &=Y^2-3XY-2Y-\frac{I}{2}e^{-6t}.\\
   \end{aligned}
  \right. 
\end{equation}

\begin{Proposition}
After a series of explicit blowing-ups at ten points including seven infinitely near points, the phase space ${\mathcal X}$ for this system can be obtained by gluing four copies of ${\mathbb C}^2 \times {\mathbb C}${\rm:\rm}
\begin{center}
${U_j} \times {\mathbb C} \cong {\mathbb C}^2 \times {\mathbb C} \ni \{(x_j,y_j,t)\},  \ \ j=0,1,2,3$
\end{center}
via the following birational transformations{\rm:\rm}
\begin{align}
\begin{split}
0) \ x_0=&x, \quad y_0=y,\\
1) \ x_1=&\frac{1}{x}, \quad y_1=\left((yx+\frac{I}{4}e^{-6t})x+\frac{I}{2}e^{-6t}\right)x,\\
2) \ x_2=&\frac{1}{x}, \quad y_2=\left(((y-2x)x-\frac{I}{4}e^{-6t})x+\frac{I}{2}e^{-6t}\right)x,\\
3) \ x_3=&xy, \quad y_3=\frac{1}{y}.\\
\end{split}
\end{align}
\end{Proposition}

\noindent
Here, for notational convenience, we have renamed $X,Y$ to $x,y$.

We remark that the phase space ${\mathcal X}$ is not a rational surface of type $E_7^{(1)}$ (see Figure 1).

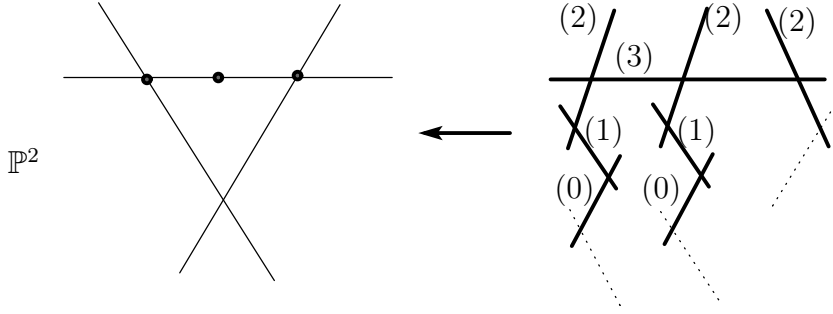
\begin{figure}
\unitlength 0.1in
\begin{picture}(42.60,17.40)(14.40,-21.30)
%
\special{pn 8}%
\special{pa 1730 780}%
\special{pa 3430 780}%
\special{fp}%
%
\special{pn 8}%
\special{pa 1910 390}%
\special{pa 2820 1840}%
\special{fp}%
%
\special{pn 8}%
\special{pa 3160 400}%
\special{pa 2330 1800}%
\special{fp}%
\put(14.4000,-12.8000){\makebox(0,0)[lb]{${\mathbb P}^2$}}%
%
\special{pn 20}%
\special{sh 0.600}%
\special{ar 2160 790 20 20  0.0000000 6.2831853}%
%
\special{pn 20}%
\special{sh 0.600}%
\special{ar 2530 780 20 20  0.0000000 6.2831853}%
%
\special{pn 20}%
\special{sh 0.600}%
\special{ar 2940 770 20 20  0.0000000 6.2831853}%
%
\special{pn 20}%
\special{pa 4040 1070}%
\special{pa 3600 1070}%
\special{fp}%
\special{sh 1}%
\special{pa 3600 1070}%
\special{pa 3667 1090}%
\special{pa 3653 1070}%
\special{pa 3667 1050}%
\special{pa 3600 1070}%
\special{fp}%
%
\special{pn 20}%
\special{pa 4250 790}%
\special{pa 5670 790}%
\special{fp}%
%
\special{pn 20}%
\special{pa 4580 430}%
\special{pa 4340 1150}%
\special{fp}%
\special{pa 4300 930}%
\special{pa 4590 1360}%
\special{fp}%
\special{pa 4610 1190}%
\special{pa 4360 1660}%
\special{fp}%
%
\special{pn 20}%
\special{pa 5060 420}%
\special{pa 4820 1140}%
\special{fp}%
\special{pa 4780 920}%
\special{pa 5070 1350}%
\special{fp}%
\special{pa 5090 1180}%
\special{pa 4840 1650}%
\special{fp}%
%
\special{pn 20}%
\special{pa 5370 430}%
\special{pa 5690 1140}%
\special{fp}%
%
\special{pn 8}%
\special{pa 4350 1470}%
\special{pa 4610 1970}%
\special{dt 0.045}%
\special{pa 4610 1970}%
\special{pa 4610 1969}%
\special{dt 0.045}%
%
\special{pn 8}%
\special{pa 4820 1480}%
\special{pa 5120 1940}%
\special{dt 0.045}%
\special{pa 5120 1940}%
\special{pa 5120 1939}%
\special{dt 0.045}%
%
\special{pn 8}%
\special{pa 5700 960}%
\special{pa 5400 1450}%
\special{dt 0.045}%
\special{pa 5400 1450}%
\special{pa 5400 1449}%
\special{dt 0.045}%
\put(45.8000,-7.6000){\makebox(0,0)[lb]{(3)}}%
\put(42.9000,-5.6000){\makebox(0,0)[lb]{(2)}}%
\put(44.2000,-11.5000){\makebox(0,0)[lb]{(1)}}%
\put(42.7000,-14.5000){\makebox(0,0)[lb]{(0)}}%
\put(50.4000,-5.6000){\makebox(0,0)[lb]{(2)}}%
\put(49.0000,-11.5000){\makebox(0,0)[lb]{(1)}}%
\put(47.2000,-14.5000){\makebox(0,0)[lb]{(0)}}%
\put(54.2000,-5.8000){\makebox(0,0)[lb]{(2)}}%
\put(41.2000,-23.0000){\makebox(0,0)[lb]{$(*)$:intersection number}}%
\end{picture}%
\label{sec1}
\caption{The bold lines denote $(-2)$-curve.}
\end{figure}

It is still an open question whether integrability status of this system is known or not.

\section{The case of $(\sigma,\varepsilon,b)=(1,-3,2)$}

\begin{equation}\label{5.1}
  \left\{
  \begin{aligned}
   \frac{dx}{dt} &=y+3x,\\
   \frac{dy}{dt} &=-xz+x+3y,\\
   \frac{dz}{dt} &=xy+6z
   \end{aligned}
  \right. 
\end{equation}

\begin{Proposition}
This system has
\begin{equation}
I:=e^{-6t}(x^2-2z)
\end{equation}
as its first integral.
\end{Proposition}


\begin{thebibliography}{99}
\bibitem[1]{1} M. J. Ablowitz and P. A. Clarkson, 
{\em Solitons, Nonlinear Evolution Equations and Inverse Scattering}, 
L. M. S. Lect. Notes Math. {\bf 149}, C. U. P., Cambridge, (1991). 


\bibitem[2]{2} J. Chazy, 
{\em Sur les \'equations diff\'erentielles du trousi\'eme ordre et d'ordre sup\'erieur dont l'int\'egrale a ses points critiques fixes}, 
Acta Math. {\bf 34} (1911), 317--385. 

\bibitem[3]{3} P. A. Clarkson and P. J. Olver, 
{\em Symmetry and the Chazy equation}, 
J. Diff. Eqns. {\bf 124}, (1996), 225--246. 

\bibitem[4]{4} C. M. Cosgrove, 
{\em Chazy classes IX-XI of third-order differential equations}, 
Stud. Appl. Math. {\bf 104}, (2000), 171--228. 

\bibitem[5]{5} C. M. Cosgrove, 
{\em Higher-order Painlev\'e equations in the polynomial class I. Bureau symbol P2}, Stud. Appl. Math. {\bf 104}, (2000), 1--65. 

\bibitem[6]{6} E. L. Ince, 
{\em Ordinary differential equations}, Dover Publications, New York, (1956). 


\bibitem[7]{7} T. Matano, A. Matumiya and K. Takano, 
{\em On some Hamiltonian structures of Painlev\'e systems, II}, 
J. Math. Soc. Japan {\bf 51} (1999), 843--866. 

\bibitem[8]{8} Y. Ohyama, {\em Differential relations of theta functions}, Osaka J. Math. {\bf 32} (1995), 431--450.  

\bibitem[9]{9} Y. Ohyama, {\em Nonlinear systems related to second order linear equations}, Osaka J. Math. {\bf 33} (1996), 927--949.  

\bibitem[10]{10} Y. Ohyama, {\em Differential equations for modular forms with level three}, Funkcial. Ekvac. {\bf 44} (2001), 377--389.  


\bibitem[11]{11} K. Okamoto, {\em Sur les 
feuilletages associ\'es aux \'equations du second ordre \`a points critiques fixes de P. Painlev\'e, Espaces des conditions initiales}, Japan. 
J. Math. {\bf 5} (1979), 1--79.  

\bibitem[12]{12} K. Okamoto, {\em Polynomial Hamiltonians associated with Painlev\'e equations}, I, II, Proc. Japan Acad. {\bf 56} (1980),  264--268; ibid, 367--371.

\bibitem[13]{13} T. Shioda and K. Takano, {\em On some Hamiltonian structures of Painlev\'e systems I}, Funkcial. Ekvac. {\bf 40} (1997), 271--291. 

\bibitem[14]{14} Y. Sasano, {\em Resolution of singularities for a family of third-order differential systems with small meromorphic solution spaces}, Kyushu J. Math. {\bf 59} (2005), 307--320. 

\bibitem[15]{15} W.-H. Steeb and N. Euler, {\em Nonlinear Evolution Equations and Painlev\'e Test}, World Scientific. (1998), ISBN 9971-50-744-7. 
\end{thebibliography}
\end{document}